\documentclass[aps, prl, twocolumn, groupedaddress, showpacs, floatfix,10pt]{revtex4-1}

\usepackage{amsmath}
\usepackage[sc]{mathpazo}
\usepackage{bm}
\usepackage{amscd}
\usepackage{amsthm}
\usepackage{amssymb} 
\usepackage{latexsym}
\usepackage{euscript}
\usepackage{epsfig}
\usepackage{graphics}
\usepackage{array} 
\usepackage{enumerate}
\usepackage{color}  
\usepackage{hyperref}    

\newcommand{\bel}[1]{\begin{equation}\label{#1}}

\newcommand{\be}{\begin{equation}}

\newcommand{\ba}{\begin{eqnarray}}
\newcommand{\ea}{\end{eqnarray}}
\newcommand{\rf}[1]{(\ref{#1})}

\newcommand{\qe}{\end{equation}}

\newcommand{\beq}{\begin{equation}}
\newcommand{\ee}{\end{equation}}
\newcommand{\benn}{\begin{equation*}}
\newcommand{\eenn}{\end{equation*}}
\newcommand{\bea}{\begin{eqnarray}}
\newcommand{\eea}{\end{eqnarray}}
\newcommand{\beann}{\begin{eqnarray*}}
\newcommand{\eeann}{\end{eqnarray*}}

\newcommand{\x}{{\mathbf x}}

\newcommand{\f}{{\mathbf f}}
\newcommand{\g}{{\mathbf g}}
\newcommand{\h}{{\mathbf h}}

\newcommand{\F}{{\mathbf F}}

\newcommand{\eps}{\epsilon}

\def\txtd{{\textnormal{d}}}
\def\txte{{\textnormal{e}}}

\def\txtD{{\textnormal{D}}}

\begin{document}

\title{Coupled Dynamics on Hypergraphs:\\
Master Stability of Steady States and Synchronization}

\author{Raffaella Mulas}
\affiliation{Max Planck Institute for Mathematics in the Sciences, Inselstr.~22, 04103 Leipzig, Germany}

\author{Christian Kuehn}
\affiliation{Faculty of Mathematics, Technical University of Munich, Boltzmannstr.~3, 85748 
Garching b.~M\"unchen, Germany\\ and Complexity Science Hub Vienna, Josefst\"adter Str.~39, 1080 Vienna, 
Austria}

\author{J\"urgen Jost} 
\affiliation{Max Planck Institute for Mathematics in the Sciences, Inselstr.~22, 04103 Leipzig, Germany\\
and Santa Fe Institute for the Sciences of Complexity, 1399 Hyde Park Road Santa Fe, New Mexico 87501, USA}

\begin{abstract}
\noindent In the study of dynamical systems on networks/graphs, a key theme  is how the network topology influences 
stability for steady states or synchronized states. Ideally, one would like to derive conditions
for stability or instability that instead of  microscopic details of the individual
nodes/vertices rather make the influence of the network coupling topology visible. The master stability function
is an important such tool to achieve this goal. Here we generalize the master
stability approach to hypergraphs. A hypergraph coupling structure is important as it allows us to 
take into account arbitrary higher-order interactions between nodes.  As for instance in the theory of coupled map lattices, we study Laplace type interaction structures in detail. Since the spectral theory of Laplacians on hypergraphs is richer than on graphs, we see the possibility of new dynamical phenomena. 
More generally, our arguments provide a blueprint for  how to generalize  dynamical structures and results
from graphs to hypergraphs.    
\end{abstract}

\maketitle

\section{Introduction}

Dynamical systems on networks are a fundamental part of the  theory of complex 
systems~\cite{BarratBarthelemyVespignani,PorterGleeson}. A common situation in network dynamics 
is that one would like to infer dynamical conclusions just from the underlying network structure.
This  has led to the introduction of the master stability function 
formalism~\cite{Pecora}, see also the exposition in~\cite{Newman}. The idea is to assume
sufficient symmetry and/or common dynamics for each individual node/vertex, which then 
makes it possible to re-write stability conditions for steady states, or even more complicated
synchronized solutions, in terms of network data. Examples of network data in this context
are spectra, e.g., of the graph Laplacian or the adjacency matrix~\cite{Chung}. The master 
stability function approach has been successfully applied in many applications, particularly
in the context of synchronization of oscillators~\cite{BarahonaPecora,ArenasDiazPerez,Nishikawaetal};
see also the surveys~\cite{Arenasetal,DoerflerBullo}.\\
However,  just considering binary interactions modelled by
a network/graph is often insufficient in applications. One then needs 
generalizations of graphs. A first natural generalization are simplicial complexes~\cite{Hatcher}.
Simplicial complexes have appeared in several applications, e.g., in protein 
classification~\cite{Cangetal}, in percolation models for statistical 
physics~\cite{BianconiKryvenZiff}, in computational neuroscience~\cite{GiustiGhristBassett}, 
 in modelling dynamics of social peer pressure~\cite{HorstmeyerKuehn}, or in epidemiology~\cite{Iacopinietal,MatamalasGomezArenas}. More generally, these results are examples
that \emph{higher-order interactions}~\cite{BickAshwinRodrigues,Grillietal,SkardalArenas} are 
relevant between nodes/vertices, where we note that the study of higher-order interactions has 
already quite a long history, particularly in ecology~\cite{Abrams,BillickCase}. While simplicial complexes form a very convenient mathematical structure, they are
also somewhat rigid as not all possible higher-order interactions are allowed.  This led to an 
interest to study more general \emph{hypergraphs}, e.g., for cellular networks~\cite{KlamtHausTheis}, 
for opinion formation~\cite{LanchierNeufer}, for epidemic spreading~\cite{BodoKatonaSimon}, or for 
social network analysis~\cite{ZhangLiu}. For instance, consider collaboration relations among scientists (see for instance \cite{PDJ}). We may have scientists $A,B,C$ that coauthor a paper, and there may also exist a paper written by $A$ and $B$ without $C$, as well as single author papers by $A$ and $C$, but no others. This would be modelled by a hypergraph with vertices $A,B,C$ and hyperedges $\{A\},\{C\},\{A,B\},\{A,B,C\}$. Neither a graph nor a simplicial complex would be adequate to capture this structure. \\
Therefore, in this paper, we study dynamics on hypergraphs. We shall  generalize the general tool of master stability functions from graphs to hypergraphs. In particular, we derive general conditions for the linear stability of synchronized dynamics.  We then turn to the important special class of Laplace type interactions, which arise in many applications, e.g., in the consensus problem~\cite{JardonKuehn1,SaberMurray}. In this context, we can apply the recently developed spectral theory for hypergraph Laplacians~\cite{JM}. At the end, we provide an outlook how our framework could be used as a 
blueprint to systematically generalize dynamical aspects of graphs to  hypergraphs. 
  
\section{Setting: Stability for Systems of ODEs}

We briefly recall linear stability theory for systems of ordinary differential 
equations to fix the notation and the main ideas. Let us consider a set of units 
$i=1,\dots, N$, called nodes or vertices, in the sequel, that are dynamically interacting 
with each other. This leads to a system of differential equations,
\bel{1}
\frac{\txtd \x_i(t)}{\txtd t}= \F_i(\x_1(t), \dots ,\x_N(t)) \text{ for } t\ge 0 ,
\qe
where we assume that the state variables $\x_i$ could be vector-valued, 
$\x_i=(x_i^1,\dots ,x_i^m)$. Hence, $\F_i$ also is a vector, 
$\F_i=(F_i^1,\dots ,F_i^m)$. We may then also write \rf{1} in matrix form
\bel{2}
\frac{\txtd \x}{\txtd t}= \F(\x).
\qe
A solution $\x^\ast$ of \rf{1} is called {\it linearly stable}, or simply {\it stable} for 
short in the sequel, if any solution $\bm{\epsilon}$ of the linearization 
\bel{4}
\frac{\txtd \bm{\epsilon}_i}{\txtd t}= \sum_{j=1}^N
\frac{\partial \F_i(\x_1^\ast, \dots ,\x_N^\ast)}{\partial \x_j}\bm{\epsilon}_j ,
\qe
or in the more abstract version corresponding to \rf{2}
\bel{4a}
\frac{\txtd \bm{\epsilon}}{\txtd t}= \txtD\F(\x^\ast) \bm{\epsilon}, 
\qe
converges to $0$ for $t\to \infty$. Here, 
$\frac{\partial \F_i(\x_1^\ast, \dots ,\x_N^\ast)}{\partial \x_j}$ 
is the vector with components 
$\frac{\partial \F_i(\x_1^\ast, \dots ,\x_N^\ast)}{\partial \x_j^\alpha}, \alpha=1,\dots ,m$, and similarly for $\bm{\epsilon}$,
and therefore, in \rf{4} there is an implicit sum over $\alpha$. Linear stability is simply 
a condition on the Lyapunov exponents of the tensor $\txtD\F(\x^\ast)$
  (note that this tensor will 
in general depend on time $t$, since we are not assuming that $\x^\ast(t)$ is constant). The 
stability condition then can be expressed in terms of a Lyapunov exponent (see for 
instance \cite{Arnold}),
\bel{6}
\limsup_{t\to \infty}\frac{1}{t}\log\|\txte^{t\txtD\F(\x^\ast(t))}\|<1.
\qe
There are two special cases that are of particular interest.
\begin{enumerate}
\item The solution $\x^\ast$ is constant in time, that is, steady or stationary. 
This means that for each $i$, $\x_i^\ast(t)=\x_i^\ast(0)$ is independent of time $t$. 
Such a stationary state simply satisfies
\bel{8} 
\F_i(\x_1^\ast, \dots ,\x_N^\ast) =0 \quad \text{ for } t\ge 0 .
\qe
For such a solution, the stability condition is simply \rf{6}. 
\item The solution $\x^\ast(t)$ represents a synchronized state. This means that it is 
independent of the vertex $i$, that is, $\x_i^\ast(t)=\x_j^\ast(t)$ for all $i$ and $j$, 
and all $t$. To make such a solution feasible, we should also assume that $\F_i$ is the 
same for all $i$. For the stability of synchronization, we only need to require that any non-synchronized  solution of \rf{4} converges to $0$ for $t\to \infty$. 
\end{enumerate}
In the sequel, we shall only consider the second case. The first case succumbs to a similar, but easier  analysis. 

\section{Interaction on Networks}
\label{net}

We now consider the situation where a vertex $i$ does not interact indiscriminately 
with all other vertices but only maintains interactions with a subset of vertices; those 
vertices are called the neighbors of $i$, and one writes $j\sim i$ when $j$ is such a neighbor 
of $i$. When one considers network interactions, these interactions are assumed to be pairwise 
only. That means that we are able to write the dynamical system \rf{1} in the form
\bel{11}
\frac{\txtd \x_i}{\txtd t}= \f_i(\x_i) +\sum_{j, j\sim i} \g_{ij}(\x_i,\x_j) \text{ for } t\ge 0 .
\qe
Here, $\f_i$ is a self-interaction term of $i$, whereas $\g_{ij}$ stands for the pairwise interaction 
between $i$ and $j$. In order to make the interaction structure more explicit, one often 
considers particular subclasses of systems of the form \rf{11} such as (see also~\cite{Newman})
\bel{12}
\frac{\txtd \x_i}{\txtd t}= \f(\x_i) +\sum_{j} a_{ij} \g(\x_i,\x_j) 
\qe
or the even simpler subclass (see also~\cite{Pecora}) 
\bel{13}
\frac{\txtd \x_i}{\txtd t}= \f(\x_i) +\sum_{j} a_{ij} \h(\x_j), 
\qe
where the (vector-valued) dynamical functions $\f ,\g ,\h$ no longer depend on the vertices. The reason to consider simpler subclasses such as \rf{12} and/or \rf{13} is 
twofold. Firstly, these structures appear frequently in modelling, e.g., in the context of 
neuroscience and for various problems regarding synchronization. Secondly, a general result 
for the stability of systems \rf{11} cannot be expected as there is too little specific 
mathematical structure, so we have to strike a balance between modelling simplifications and
obtainable theoretical results. The forms \rf{12} and \rf{13} have shown to be very useful in
the context of graphs~\cite{Newman,Pecora}, so they form a natural starting point for an
extension to hypergraphs. Based on these considerations, the focus then is on the interaction 
matrix $A=(a_{ij})_{i,j=1,\dots N}$. The neighborhood 
structure can be included in that matrix by stipulating that $a_{ij}=0$ unless $j\sim i$. 

We consider \rf{13}, as the analysis of \rf{12} is similar. The resulting  stability condition has been referred to in the literature as master stability condition. 
If one wishes to make synchronized dynamics possible, one usually assumes that
\bel{14}
a:=\sum_{j} a_{ij}
\qe
does not depend on $i$. In that case, a synchronized solution $\x^\ast$ of \rf{13} would satisfy
\bel{15}
\frac{\txtd \x^\ast(t)}{\txtd t}= \f(\x^\ast(t)) + a \h(\x^\ast(t)).  
\qe
The linear stability equation \rf{4a} for \rf{13} at a solution $\x^\ast$ is (\cite{Pecora})
\bel{16}
\frac{\txtd \bm{\epsilon}}{\txtd t}= (\textnormal{Id}\otimes \txtD\f(\x^\ast) +A\otimes \txtD\h(\x^\ast))\bm{\epsilon} ,
\qe
where $\textnormal{Id}$ always denotes the identity operator of suitable size, which is simply the 
$N$-dimensional identity matrix in the context of \rf{16}. 
When we assume that the coupling matrix $A$  can be diagonalized (for instance, if it is 
symmetric, i.e., $a_{ij}=a_{ji}$ for all $i,j$), we let its eigenvalues be $\mu_k, k=1,\dots ,N$. 
Since $\textnormal{Id}$ is the identity matrix, we can decompose \rf{16} into the corresponding modes 
$\eps^k$, that is, 
\bel{18}
\frac{\txtd \eps^k}{\txtd t}= (\txtD\f(\x^\ast) +\mu_k  \txtD\h(\x^\ast))\eps^k  
\text{ for } k=1,\dots ,N.
\qe
When we assume \rf{14}, one of the eigenvectors of $A$ is constant. Therefore, at a synchronized 
state $\x^\ast$, we obtain a mode $\eps^1(t)$ with $\eps_i^1(t)=\eps_j^1(t)$ for all $i,j$. The 
evolution of the mode therefore leaves the synchronization manifold invariant. 
Synchronization is stable when all other modes decay. Let us consider the case where $\f=\h$. 
Then \rf{18} becomes
\bel{19}
\frac{\txtd \eps^k}{\txtd t}= (1+\mu_k) \txtD\f(\x^\ast)\eps^k  \text{ for } k=1,\dots ,N.
\qe
The stability condition then is (see \cite{JJ1})
\bel{20}
\limsup_{t\to \infty}\frac{1}{t}\log \|\txte^{t(1+\mu_k) \txtD\f(\x^\ast)}\| <1,
\qe
that is,
\bel{21}
|1+\mu_k| \ell_{\f} <1 
\qe
where as in \rf{6}, 
\bel{22}
\ell_{\f}:= \limsup_{t\to \infty}\frac{1}{t}\log \|\txte^{t \txtD\f(\x^\ast)}\|
\qe
is the maximal Lyapunov exponent of $\f$ (at the particular solution $\x^\ast$, but in order 
to have a general criterion, we may take the supremum over all solutions). The inequality \rf{21} now 
separates and relates the condition for the dynamical update $\f$ and the network connectivity 
as encoded in the coupling matrix $A$ and its eigenvalues. In the interesting case, we have 
$\ell_{\f} >1$, that is, the dynamics generated by $\f$ is unstable. But if the eigenvalues 
$\mu_2,\dots ,\mu_N$ lie between $-2$ and $0$ and satisfy \rf{21}, synchronization may still 
be a stable state. Similar to \cite{JJ1}, we now consider the case where
\bel{23}
\frac{\txtd \x_i}{\txtd t}=  \f(\x_i) - \sigma (\Delta \f)(\x_i). 
\qe
Here, $0\le \sigma \le 1$ is a parameter and
\bel{24}
(\Delta u)(\x_i):=u(\x_i)- \frac{1}{\deg i}\sum_{j\sim i} u(\x_j) 
\qe
is the normalized Laplace operator of the network (see for instance \cite{Chung,J1} for the 
theory, but note that the conventions employed here are somewhat different from those in these 
references). The eigenvalues of $\Delta$ satisfy
\bel{26}
0=\lambda_1 \le \lambda_2 \le \dots \le \lambda_N \le 2,
\qe
where the eigenfunction for $\lambda_1=0$ is constant. The stability condition \rf{21} then 
becomes
\bel{27}
|1-\sigma \lambda_k|\ell_{\f} <1 \text{ for } k=2,\dots ,N,
\qe
that is, by \rf{26},
\bel{28} 
\lambda_2> \frac{1-\ell_{\f}^{-1}}{\sigma} \text{ and } 
\lambda_N< \frac{1+\ell_{\f}^{-1}}{\sigma}.
\qe
Thus, we need at the same time a lower bound for the first nonzero eigenvalue and an 
upper bound for the largest eigenvalue. $\lambda_2$ is controlled from below by the 
so-called Cheeger inequality \cite{Alon,Dodziuk} which quantifies the cohesion of the 
graph. $\lambda_2$ is largest when the graph is complete, and of course, a complete 
graph is more conducive to synchronized dynamics than a less coherent one. In particular, 
$\lambda_2=0$ precisely if the graph is disconnected, and for such a graph, we obviously 
cannot expect dynamics to synchronize. In fact, when the graph has more than one component, the dynamics could be synchronized on each component, but not necessarily between components. Let us consider the case of two components $\Gamma_1,\Gamma_2$. An eigenfunction for $\lambda_2=0$ then is constant on each component (with the weighted sum of the constants being zero). When $\ell_\f >1$, but \rf{28} is satisfied now for $\lambda_3$, then what we may call the generalized synchronization manifold, that is, the family of dynamical states that are synchronized inside the two components only, is stable against perturbations by other eigenstates. Analogously, of course, for more than two components.-- $\lambda_N=2$ holds precisely if the graph is bipartite, 
and in fact the gap $2-\lambda_N$ quantifies the deviation from bipartiteness \cite{BJ}. 
On a bipartite graph, antiphase oscillations are possible, and thus, there again is an 
obstacle to synchronization carried by the mode associated with $\lambda_N$. That is why we 
need the upper bound. Given $\ell_{\f}$ and the topology of the underlying graph, \rf{28} 
then tells us whether we can find a range of coupling strengths $\sigma$ for which 
synchronized dynamics are stable. 

\section{Interaction on hypergraphs}
\label{hyp}
So far, we have essentially summarized or reformulated known results. 
In particular, in the preceding section, we have considered dynamics on a network where 
the dynamics at each vertex is coupled with the dynamics of its neighbors. The network 
thus corresponds to a graph with edges defined by the neighborhood relations. Thus, all 
relations are binary. When we also want to include higher order interactions, as in many 
empirical systems, we need an underlying structure that is more general than that of a 
graph. We need a {\it hypergraph}. A hypergraph has a set $V$ of vertices $i=1,\dots ,N$ 
and a set $H \subset 2^V$ of hyperedges $h=1, \dots , M$. Thus, each hyperedge is a set 
of vertices $h=\{ i_{h(1)}, \dots ,i_{h(m_h)}\}$ where $m_h$ is the number of vertices contained 
in the hyperedge $h$. We can then consider types of dynamics analogous to those in 
equations \rf{11}. These can  be written as
\bel{31}
\frac{\txtd \x_i}{\txtd t}= \f(\x_i) +\sum_{h: i\in h}  \g_{ih}(\x_{i_{h(1)}}, 
\dots ,  \x_{i_{h(m_h)}}).
\qe
We note that the number of arguments of an interaction function $\g_{ih}$ now depends on the 
size $m_h$ of the hyperedge $h$. When we linearize \rf{31}, we therefore need  the $N\times M$ incidence 
matrix $\mathcal{I}:=(\mathcal{I}_{ih})$ defined  by 
\begin{equation*} 
\mathcal{I}_{ih}:=\begin{cases} 1 & \text{if }  i\in h\\ 0 & \text{otherwise.} \end{cases} 
\end{equation*} 
We observe that, for each $i$ and $j$,
\begin{equation*}
\mathcal{I}_{ih}\cdot \mathcal{I}_{jh}=\begin{cases} 1 & \text{if }  
i, j\in h\\ 0 & \text{otherwise.} \end{cases}\end{equation*}Therefore,
\begin{equation}\label{31b} 
\bigl(\mathcal{I}\cdot\mathcal{I}^\top\bigr)_{ij}=\sum_{h=1}^M
\mathcal{I}_{ih}\cdot \mathcal{I}^\top_{hj}=\sum_{h=1}^M\mathcal{I}_{ih}
\cdot \mathcal{I}_{jh}=\bigl| h: i, j\in h \bigr|.
\end{equation}
Returning to the general system \rf{31}, its linearized version  at a solution 
$\x^\ast$ then is
\ba
\nonumber
\frac{\txtd \bm{\epsilon}_i}{\txtd t}&= \frac{\partial \f(\x^\ast_i)}{\partial \x_i} 
\bm{\epsilon}_i +\sum_{h: i\in h} \sum_{j\in h}\frac{\partial \g_{ih}(\x^\ast_{i_{h(1)}}, 
\dots ,  \x^\ast_{i_{h(m_h)}})}{\partial \x_j} \bm{\epsilon}_j\\
\label{33}
=&\frac{\partial \f(\x^\ast_i)}{\partial \x_i} 
\bm{\epsilon}_i + \sum_{h: i\in h}\sum_{j}\mathcal{I}_{jh}\frac{\partial \g_{ih}(\x^\ast_{i_{h(1)}}, 
\dots ,  \x^\ast_{i_{h(m_h)}})}{\partial \x_j} \bm{\epsilon}_j.
\ea
For the stability of $\x^\ast$, we need to check as before whether $\bm{\epsilon}(t) \to 0$ 
as $t\to \infty$ for any solution of \rf{33}. 

After this general result, we now want to discuss the possibility and the 
stability of synchronized dynamics on hypergraphs. 
When we want to consider the analogue of \rf{12} or \rf{13} and again assume uniform interaction functions, these functions will now still depend on the size of the hyperedegs, as the number of their 
arguments varies with the size $m$ of the underlying hyperedge. Thus, we have functions $\g_m$.  When we have an interaction matrix $A=a_{ih}$, the dynamics then are of the form
\bel{31a}
\frac{\txtd \x_i}{\txtd t}= \f(\x_i) +\sum_{h: i\in h} a_{ih} \g_{m_h}(\x_{i_{h(1)}}, 
\dots ,  \x_{i_{h(m_h)}}).
\qe
When, for instance $a_{ih}=\mathcal{I}_{ih}$, \rf{33} becomes
\bel{32a}
\frac{\txtd \bm{\epsilon}_i}{\txtd t}=\frac{\partial \f(\x^\ast_i)}{\partial \x_i} 
\bm{\epsilon}_i + \sum_{j, h}\mathcal{I}_{ih} \mathcal{I}_{jh}\frac{\partial \g_{m_h}(\x^\ast_{i_{h(1)}}, 
\dots ,  \x^\ast_{i_{h(m_h)}})}{\partial \x_j} \bm{\epsilon}_j.
\qe
We thus see \eqref{31b} in action. Furthermore, we require the analogue of \rf{14}, that is, $a:=\sum_h \mathcal{I}_{ih}$ does not depend on $i$.\newline

As explained already for the case of graphs, it is necessary to 
make additional assumptions to obtain a theoretically tractable, yet interesting
and applicable coupling structure. Hence, we consider the case where $\g_m(y_1,\dots ,y_m)$ 
is a normalized symmetric function 
of its entries, for instance
\begin{equation*}
 \g_m(y_1,\dots ,y_m)=\g\left(\frac{1}{m}\sum_{j=1}^m y_j\right) 
\end{equation*}
or
\begin{equation*}
 \g_m(y_1,\dots ,y_m)=\g\left((\prod_{j=1}^m y_j)^{1/m}\right) 
\end{equation*}
for some function $\g$ (when the entries are vectors, as considered here, 
these functions can be evaluated component-wise).


Importantly, we can again consider a Laplacian type coupling. The corresponding hypergraph Laplacian was constructed in \cite{JM}, where the authors worked on the more general setting of {\it chemical hypergraphs}. Here we choose to work on this more general setting, as this offers more possibilities of modelling, and we recall some properties of the corresponding Laplacian. A chemical hypergraph is given by a collection of vertices $i=1,\dots ,N$ and a collection of {\it oriented hyperedges} $h=1,\dots ,M$.
An oriented hyperedge is a non-empty ordered subset 
$(V_h,W_h)$ of $2^V\times 2^V$. The vertices in $V_h$ and $W_h$ are called the inputs 
and outputs of $h$. Changing the orientation of $h$ simply means replacing $(V_h,W_h)$ 
by $(W_h,V_h)$. $V_h$ and $W_h$ need not be disjoint, and the vertices in $V_h\cap W_h$ 
are called {\it catalysts} of $h$. The hypergraph Laplacian of \cite{JM} then is defined as
$\tilde{\Delta} u(\x_i):=$
\beann 
&\frac{\sum_{h_{\text{in}}: i\text{ input}}\biggl(\sum_{i' 
\text{ input of }h_{\text{in}}}u(\x_{i'})-\sum_{j' \text{ output of }h_{\text{in}}}
u(\x_{j'})\biggr)}{\deg i}+\\
&-\frac{\sum_{h_{\text{out}}: i\text{ output}}\biggl(\sum_{\hat{i} \text{ input of }
h_{\text{out}}}u(\x_{\hat{i}})-\sum_{\hat{j} \text{ output of }h_{\text{out}}}
u(\x_{\hat{j}})\biggr)}{\deg i}.
\eeann
This definition is invariant under changes of orientation of hyperedges. 
For a graph, an oriented edge is simply a pair of vertices, and the definition of
the hypergraph Laplacian reduces to \rf{24}. Also, chemical hypergraphs that have either only inputs or only outputs correspond to classical hypergraphs with no orientation.\\
As before, the stability condition
couples the Lyapunov exponent of the dynamical nonlinearity $\f$, the structure of the 
hypergraph as encoded by the eigenvalues $\tilde{\lambda}_k$ of $\tilde{\Delta}$, and the 
coupling parameter $\sigma$. Indeed, if we replace in \rf{23} the usual graph Laplacian 
by the hypergraph Laplacian $\tilde{\Delta}$, then we get a stability condition
\beq
\label{eq:stabHG}
|1-\sigma \tilde{\lambda}_k|\ell_{\f} <1 \text{ for } k=1,\dots ,N,
\qe
Note carefully, that although we have 
\bel{26a}
0\leq \tilde{\lambda}_1 \le \tilde{\lambda}_2 \le \dots \le \tilde{\lambda}_N,
\qe
we do not have the same strong bounds as for the usual graph Laplacian as presented 
in~\rf{26}. Yet, we can still re-write \rf{26a} as,
\bel{28a} 
\tilde{\lambda}_{min}> \frac{1-\ell_{\f}^{-1}}{\sigma} \text{ and } 
\tilde{\lambda}_N< \frac{1+\ell_{\f}^{-1}}{\sigma},
\qe
where $\tilde{\lambda}_{min}$ is the smallest non-zero eigenvalue. Even for a connected hypergraph,  $\tilde{\lambda}_2$ need not be 
greater than $0$. This, in fact, leads to an interesting class of dynamics. Let us assume that $\tilde{\lambda}_1=\dots =\tilde{\lambda}_k=0$, but $\tilde{\lambda}_{k+1}$ satisfies \rf{28a}, that is, $\tilde{\lambda}_{k+1}> \frac{1-\ell_{\f}^{-1}}{\sigma}$. Then the class of dynamics that belong to eigenstates of the Laplacian for the eigenvalue $\tilde{\lambda}=0$ is stable. This class can be larger than the locally synchronized dynamics. For instance, consider a graph with three vertices $1,2,3$ and a single hyperedge with $V=\{1\}, W=\{2,3\}$. One eigenstate for $\tilde{\lambda}=0$ is constant, but another one is given by $u(1)=1, u(2)=u(3)=\frac{1}{2}$. This would correspond to a dynamical state $\x^\ast$ with $\g(\x_2^\ast)=\g(\x_3^\ast)=\frac{1}{2}\g(\x_1^\ast)$, which would be stable under our conditions. That is, the dynamical activity at $1$ is equally split into the activities at $2$ and $3$, as prescribed by the topology of the hypergraph. -- Conversely, it may also happen that all eigenvalues of a hypergraph are positive. Take, for instance, again three vertices, and for each $i$ a hyperedge $h_i$ with $V_{h_i}=\{i\}, W_{h_i}=\{i+1,i+2\}$, counting the vertices mod~$3$. Then all eigenvalues are positive, see \cite{JM}, precluding the possibility of  synchronized dynamics.  Furthermore, another difference with the graph case 
is that $2$ does not give an upper bound to $\tilde{\lambda}_N$. In fact, $\tilde{\lambda}_N$ 
is equal to $N$ in some cases and it is not known yet whether this is the largest possible 
value for $\lambda_N$. Nevertheless, the geometrical meaning of the largest eigenvalue does not 
change. It is in fact known that, given a hypergraph $\Gamma$ with largest eigenvalue 
$\tilde{\lambda}_N$, then
\begin{equation*}
\tilde{\lambda}_N\leq\tilde{\lambda}'_N,
\end{equation*}
where $\tilde{\lambda}'_N$ is the largest eigenvalue of a \emph{bipartite hypergraph} that has 
the same number of hyperedges as $\Gamma$ and also the same number of inputs and the same number 
of outputs in each hyperedge (catalysts are not included). Also, the equality holds if and only 
if $\Gamma$ is bipartite. 

In summary, we find that once the hypergraph Laplacian appears in the dynamics directly, one
can still derive a master stability condition. But one has to
be careful, e.g., in treating the dimension of the synchronization manifold as well as 
possible degenerate additional neutral modes associated to zero eigenvalues, which may 
appear on a linear level for the hypergraph Laplacian. In addition, it is clear that
hypergraph coupling can shift the stability regions. This lends some interest to  results for a
particular model in the special case of simplicial complexes~\cite{Iacopinietal}. However, note that our master stability conditions only operate 
on the level of the linearization. The case of higher-order interactions and bifurcations,
where nonlinearities matter even locally, is far more involved~\cite{KuehnBick}.\\
Finally, we point out that while in this article, we have considered time-continuous dynamics, our scheme also applies to time-discrete dynamics. For instance, one can study the phenomenon of the synchronization of chaos \cite{Kaneko1984} on analogues of coupled map lattices on hypergraphs.

\section{Conclusion and Outlook}

In this work we have shown how to extend the master stability function framework from graphs 
to hypergraphs. In particular, we noticed how the spectral properties of the hypergraph Laplacian
enter the stability condition, and how this changes the statements we may make regarding the
interplay between network topology and dynamics. For example, it is now possible that the upper
bound on the largest eigenvalue grows significantly, while already the smallest eigenvalue can be 
bigger than zero. Conversely, even for connected hypergraphs, the multiplicity of the eigenvalue $0$ can be larger than $1$, and this leads to interesting new classes of dynamics that are more general than synchronization, but may still be locally stable under appropriate conditions. Furthermore, we found that the incidence matrix plays an important role in 
hypergraph dynamics, and it interacts in a non-trivial way with the master stability condition(s).\\ 
We point out that the approach we have taken here provides a general strategy for lifting results about dynamics on graphs to hypergraphs. The key is to identify the steps
where the adjacency matrix or the graph Laplacian play key roles, and then replace them with
analogous hypergraph objects. The spectral theory of hypergraphs is richer than that of graphs, and that lead us to identify new classes of dynamics that are more general than synchronization but for which we can still derive stability conditions analogous to those for synchronized dynamics on graphs.\\ 

{\bf Acknowledgments.}
The authors are grateful to the anonymous referees for the constructive comments. CK acknowledges support via a Lichtenberg Professorship as well as support 
via the TiPES project funded the European Union’s Horizon 2020 research and 
innovation programme under grant agreement No.~820970.

	\bibliographystyle{unsrt}
	\bibliography{Master20-06-03}

\end{document}